\documentclass{amsart}

\usepackage{amsmath,amssymb,amsthm,graphicx,slashed,subcaption}%
\usepackage{amsmath}%
\usepackage{amsfonts}%
\usepackage{amssymb}%
\usepackage{float}
\usepackage{tikz,pgfplots}
\pgfplotsset{compat=1.11}

\usetikzlibrary{patterns}
\usepackage{mathpple}
\newcommand{\rvline}{\hspace*{-\arraycolsep}\vline\hspace*{-\arraycolsep}}

\usepackage{hyperref,enumerate}

\usepackage{colortbl}

\usepackage[all,cmtip]{xy}

\definecolor{Red}{rgb}{0.6,0,0}

\usepackage{graphicx}
\providecommand{\U}[1]{\protect\rule{.1in}{.1in}}
\newtheorem{thm}{Theorem}[]
\newtheorem{corl}[thm]{Corollary}

\newtheorem{lma}[thm]{Lemma}

\newtheorem{prop}[thm]{Proposition}
\newtheorem{defn}[thm]{Definition}
\newtheorem{ex}[thm]{Example}
\newtheorem{rem}[thm]{Remark}
\def\tilde{\widetilde}

\def\C{\mathbb{C}}

\DeclareMathOperator{\dom}{Dom}
\def\E{\mathcal{E}}
\def\env{\mathrm{env}}
\def\epsilon{\varepsilon}
\def\eps{\varepsilon}

\def\H{\mathcal{H}}

\DeclareMathOperator{\ind}{Ind}
\def\I{\mathbb{I}}
\def\id{\textup{id}}
\def\K{\mathcal{K}}

\def\min{\mathrm{min}}

\def\nn{\nonumber}

\DeclareMathOperator{\Sig}{Sig}

\def\T{\mathbb{T}}

\newcommand{\Toep}[1]
           {\textup{Toep}_{#1}(\C)}
\def\U{\mathcal{U}}
\def\V{\mathcal{V}}
\def\Z{\mathbb{Z}}

\newtheoremstyle{commentstyle}
  {0.2cm}{0.2cm}
  {\sf}
  {0cm}
  {\bfseries}{ }
  {0cm}
  {\thmname{#1}\thmnumber{ #2}:\thmnote{ #3}}

\theoremstyle{commentstyle}
\newtheorem{mycomment}{Comment}

\title{A generalization of K-theory to operator systems}

\author{Walter D. van Suijlekom}

\address{Institute for Mathematics, Astrophysics and Particle Physics, Radboud
University Nijmegen, Heyendaalseweg 135, 6525 AJ Nijmegen, The Netherlands.
}
\email{waltervs@math.ru.nl}
\date{September 4, 2024}

\begin{document}
\maketitle

\begin{abstract}
  We propose a generalization of K-theory to operator systems. Motivated by spectral truncations of noncommutative spaces described by $C^*$-algebras and inspired by the realization of the K-theory of a $C^*$-algebra as the Witt group of hermitian forms, we introduce new operator system invariants indexed by the corresponding matrix size. A direct system is constructed whose direct limit possesses a semigroup structure, and we define the $K_0$-group as the corresponding Grothendieck group. This is an invariant of unital operator systems, and, more generally, an invariant up to Morita equivalence of operator systems. For $C^*$-algebras it reduces to the usual definition. We illustrate our invariant by means of the spectral localizer.  
  \end{abstract}
\tableofcontents

\section{Introduction}
The last few years has seen many new interactions between noncommutative geometry and operator theory. 
In particular, the development of noncommutative geometry to describe spectral truncations \cite{ALM14,CS20}, the corresponding spectral localizer \cite{LS18a,LS18b} ({\em cf.} \cite{DSW23} and references therein), and spaces up to tolerance relations \cite{CS21,GS22,ALL22} has led to the development of many interesting new structures and applications in operator system theory. These include operator system duality \cite{Far21,Ng21,FB23}, noncommutative convex geometry \cite{KKM21}, quantum metric spaces \cite{Sui21,Hek21,LS23,GS22,Rie23} and the notion of Morita equivalence for operator systems \cite{EKT21}.

A key ingredient for many applications in noncommutative geometry based on $C^*$-algebras is K-theory. Given the above fruitful interactions based on the replacement of $C^*$-algebras by operator systems, one is naturally led to the question whether there is an analogue of K-theory for operator systems. We address this question in the present paper.

On our wish list for a candidate for K-theory we have put that
\begin{enumerate}
\item it should capture the spectral localizer \cite{LS18a,LS18b} for spectral truncations as a pairing with K-theory;
\item it should capture compressions of the element $Y$ describing quanta of geometry in \cite{CCM14};
\item it should be invariant under Morita equivalence \cite{EKT21};
\item it should be a refinement of K-theory for $C^*$-algebras, taking the operator system structure into account.
  \end{enumerate}
Some results in the direction of K-theory for operator systems already exist in the literature. For instance, projections in operator systems were defined in \cite{AR22}. However, it is not difficult to see that spectral compressions $PpP$ of a projection $p$ in a $C^*$-algebra are in general not projections. As a matter of fact, in general they are not even $\epsilon$-projections in the sense of \cite{OY15} so that also quantitative K-theory may not be suitable for these kind of applications (instead, quantitative K-theory is more fitting for the examples coming from tolerance relations). The K-theory for absolute matrix order unit spaces of \cite{KK21} need an additional structure (the absolute value map) on the operator system which we would like to avoid.

The aforementioned applications based on the spectral localizer suggests that one should look for invertible, or non-degenerate elements instead.
In fact, this is nicely aligned with the elegant description of $C^*$-algebraic K-theory $K_0(A)$ as the Witt group of hermitian forms over $A$ (see \cite{Knu91,Ros95}). Indeed, these are given by invertible self-adjoint elements in $M_n(A)$. It is precisely the matrix structure of the operator system that allows us to formulate such a notion of hermitian forms for operator systems as well. Based on this, we will develop  and analyze K-theoretic invariants for operator sysetems.

\bigskip

This paper is structured as follows. After some brief background on operator systems and unital completely positive (ucp) maps, we introduce in Section \ref{sect:K1} a notion of hermitian forms in a unital operator system $E$. These are given by non-degenerate self-adjoint elements $x \in M_n(E)$.
A natural notion of homotopy equivalence allows to introduce sets $\V(E,n)$ of equivalence classes up to homotopy, which are labeled by the matrix size $n$. We will show that they are invariants for the operator system structure and compatible with the notion of direct sum. A semigroup of hermitian forms then appears when we introduce maps $\V(E,n) \to \V(E,n+1)$ and consider the direct limit $\varinjlim \V(E,n)$. The $K_0$-group of the operator system is then defined to be the Grothendieck group of this semigroup. We show that it coincides with $C^*$-algebraic K-theory in case the unital operator system is a unital $C^*$-algebra, and that it is stably equivalent when taking matrix amplifications. It is not to be expected that this notion of $K_0$-theory behaves nicely with respect to ucp maps, as can be seen already in the simple case of a ucp map given by compression by a projection. Nevertheless, for some maps ---including complete order isomorphisms--- we have induced maps between the corresponding $K_0$-groups.

In Section \ref{sect:non-unit} we consider the extension of $K$-theory to non-unital operator systems in the sense of \cite{Wer02}. This allows to formulate and prove the stability of $K_0$-theory, which as a consequence of \cite{EKT21} then shows invariance of $K_0$ under Morita equivalence. 
In Section \ref{sect:appl} we return to our initial motivation and formulate the spectral localizer \cite{LS18a} in terms of the K-theory groups for the spectral truncations. An illustrative example is given by spectral compressions of projections on the torus.

\subsection*{Acknowledgements}
I would like to thank Alain Connes, Jens Kaad, Matthew Kennedy, Marc Rieffel, and Steffen Sagave for fruitful discussions and suggestions. I thank Malte Leimbach for a very careful proofreading of a draft of this paper.

\section{Background on operator systems}

We start by briefly recalling the theory of operator systems, referring to \cite{ER00,Pau02,Pis03,Bla06} for more details.

A {\em unital operator system} $(E,e)$ is a matrix-ordered $*$-vector space $E$, equipped with an Archimedean order unit $e$. A map $\phi: E \to F$ between operator systems determines a family of maps $\phi^{(n)} : M_n(E) \to M_n(F)$ given by $\phi^{(n)}([x_{ij}]) = [\phi(x_{ij})]$. A map $\phi: E \to F$ between unital operator systems is called {\em completely positive} if each $\phi^{(n)}$ is positive ($n \geq 1$). We also abbreviate  completely positive by {\em cp}, and unital completely positive by {\em ucp}. 

A {\em dilation} of a ucp map $\phi: E \to B(\H)$ of a unital operator system is a ucp map $\psi: E \to B(\K)$, where $\K$ is a Hilbert space containing $\H$ such that $P_\H \psi(x)|_\H = \phi(x)$ for all $x \in E$. 
The ucp map $\phi$ is called {\em maximal} if every dilation of $\phi$ is obtained by attaching a direct summand.

A non-zero cp map $\phi:E \to B(\H)$ is said to be {\em pure} if the only cp maps satisfying $0 \leq \psi \leq \phi$ are scalar multiples of $\phi$.

We may view $E$ as a concrete operator system in the $C^*$-algebra $C^*(E)$ it generates; in this case, we say that a ucp map $\phi: E \to B(\H)$ has the {\em unique extension property} if it has a unique ucp extension to $C^*(E)$ which is a $*$-representation. If, in addition, the $*$-representation is irreducible, it is called a {\em boundary representation} \cite{Arv69}. The following result is well-known in the literature \cite{Arv69,MS98,Arv08,Kle14,DK15} but for completeness we include a proof. 
\begin{prop}
  \label{prop:bdry-pure-max}
  Let $\phi: E \to B(\H)$ be a ucp map. Then
  \begin{enumerate}
  \item $\phi $ is maximal if and only if it has the unique extension property.
    \item $\phi$ is pure and maximal if and only if it is a boundary representation. 
    \end{enumerate}
\end{prop}
\proof
(1) is \cite[Proposition 2.4]{Arv08}.
(2) If the extension $\tilde \phi: C^*(E) \to B(\H)$ of $\phi$ is reducible, then there exists a non-trivial projection $P \in \tilde \phi(C^*(E))' $. Then the map $\psi:E \to B(\H)$ defined by $\psi(x) = P \tilde \phi(x)$ is cp and $0 \leq \psi \leq \phi$. But $\psi(1) = P$ while $\phi(1) = \id_\H$ so that $\psi$ is not a scalar multiple of $\phi$; hence $\phi$ is not pure.

For the other implication, suppose that $\phi$ is a boundary representation and take a cp map $\psi$ such that $0 \leq \psi \leq \phi$. By Arveson's extension theorem, we have cp maps $\tilde \psi$ and $\widetilde{\phi-\psi}$ from $C^*(E)$ to $B(\H)$. Since $\tilde \psi + \widetilde{\phi-\psi}$ extends $\phi$ it follows from the unique extension property of $\phi$ that $\widetilde{\phi-\psi} =\tilde \phi - \tilde \psi$. Consequently, $\tilde \psi \leq \tilde\phi$ so that by \cite[Theorem 1.4.2]{Arv69} there exists an operator $T \in \tilde \phi(C^*(E))'$ such that
$$
\tilde \psi(a) = T \tilde \phi(a) ; \qquad \forall a\in C^*(E).
$$
Since $\tilde \phi$ is irreducible, it follows that $T = t \cdot \id_\H$ for some $t \in [0,1]$ so that $\tilde \psi  = t \tilde \phi$ and, consequently, also $\psi = t \phi$. Hence $\phi$ is pure.      
\endproof

\section{K-theory for unital operator systems}
\label{sect:K1}
Even though in an operator system we cannot speak about invertible elements, we may use the pure and maximal ucp maps to introduce the following notion of nondegeneracy. 
\begin{defn}
\label{defn:herm-forms}
Let $(E,e)$ be a unital operator system.  A self-adjoint element $x=x^* \in M_n(E)$ is called a {\em hermitian form} if it is non-degenerate in the sense that there exists $g>0$ such that for all pure and maximal ucp maps  $\phi: E \to B(\H)$ we have
\begin{equation}
|\phi^{(n)}(x)| \geq g \cdot \id_\H^{\oplus n}
  \label{eq:nondeg}
\end{equation}
The smallest real number $g>0$ such that \eqref{eq:nondeg} holds is called the {\em gap} of $x$.
\end{defn}

We denote the set of hermitian forms contained in $M_n(E)$ by $H(E,n)$
. 
The following result will be of crucial importance to us, as it makes it feasible to check the non-degeneracy condition in concrete cases:
\begin{prop}
A self-adjoint element $x \in M_n(E)$ is a hermitian form if and only if $\imath_E^{(n)}(x)$ is an invertible element in the $C^*$-envelope $C^*_\env(E)$. Moreover, $x$ has gap $g>0$ if and only if $|\imath_E^{(n)}(x)| \geq g \cdot 1^{\oplus n}_{C^*_\env(E)}$.
  \end{prop}
\proof
In \cite{DK15} the $C^*$-envelope of $E$ is constructed as the direct sum of all boundary representations $(\H_\sigma, \sigma)$:
$$
\imath_E : E \to \bigoplus_\sigma B( \H_\sigma ).
$$
Now if $x \in E$ then $\imath_E(x)$ is invertible if and only if $| \imath_E(x)|\geq g \cdot \id_{\oplus_\sigma (\H_\sigma)}$, which holds if and only if $|\sigma(x)| \geq g \cdot \id_{\H_\sigma}$ for all boundary representations $\sigma$. A similar statement holds if $x \in M_n(E)$. Since by Proposition \ref{prop:bdry-pure-max} a ucp map is a boundary representation if and only if it is pure and maximal, the result follows.
\endproof

A trivial example of a hermitian form is the order unit $x = e$ (with gap $g =1$). Other examples that in fact motivated the above definition are

\begin{ex}
\begin{enumerate}[(i)]
\item (Hermitian forms and Witt groups)  The relation between the K-theory and the Witt group for rings \cite{Knu91,Bal05} and (unital) $C^*$-algebras \cite{Ros95} stresses the role played by hermitian forms. In fact, any hermitian form on a finitely generated projective module of the form $e A^n$ is described by a self-adjoint element $x \in e M_n(A)e$. The fact that the form is non-degenerate translates to invertibility of the hermitian element $h = x+ (1-e) \in M_n(A)$, so that $h^2 \geq g^2 \cdot 1_A$ for some $g>0$. But then 
  \begin{align}\nn
    x^2  & = (ehe)^2 = eh^2 e  \geq g^2 \cdot e,
    \label{eq:trick}
  \end{align}
  since $h$ and $e$ commute. We find that $x$ is a hermitian form in the operator system $(eM_n(A)e,e)$  with gap $g$ (note that since $eM_n(A)e$ is a $C^*$-algebra, it coincides with its $C^*$-envelope). This should also explain the above terminology.

Note that in this case one may just as well consider the invertible element $x+ (1-e)$ as a hermitian form in $M_n(A)$.

\item (Almost projections and quantitative $K$-theory) In \cite{OY15} quantitative (aka controlled) K-theory $K_0^{\epsilon,r}(A)$ of a filtered $C^*$-algebra $A = (A_r)$ was defined in terms of $\epsilon$-$r$-projections, {\em i.e.} elements $p \in A_r$ such that $\| p^2 - p\| < \epsilon$ where $\epsilon < 1/4$. But then $x= 1_A -2p$ is a hermitian form in the unital operator system $(A_r,1_A)$ with gap $g = \sqrt{1-4 \epsilon} >0$ since
  $$
x^2 = 1_A - 4 (p-p^2) \geq (1-4 \epsilon )\cdot 1_A
$$

\item (Projections in operator systems)
  In \cite{AR22} projections in unital operator systems are abstractly defined. These then turn out to be precisely those elements $p \in E$ such that $\imath_E(p)$ is a projection in the $C^*$-envelope. If we set $x= e-2p$ we find that $\imath_E(x)^2 =e$ so that again such projections define hermitian forms (with gap $1$).

\item (The even spectral localizer) The even spectral localizer \cite{LS18a} is defined as a spectral compression $x = P H P$ of an invertible self-adjoint element $H$ in a $C^*$-algebra $A \subseteq B(\H)$ by a projection in $\H$. If the spectrum of $H$ does not intersect with the interval $[-g,g]$ then we compute for the operator system $(PAP,P)$ that
  \begin{equation}
  \label{eq:trick} 
x^2 =  P H P P H P =P H^2P + P H[P,H] P  = P H^2P + P [P,H][P,H] P 
  \end{equation}
  Hence if we set $\delta =\| [P,H]\|$ we find that  $x^2 \geq (g^2 -\delta^2) P$ in the $C^*$-extension $C^*(PAP) \subseteq B(\H)$. By the universal property of the $C^*$-envelope, there is a surjective $*$-homomorphism $\rho: C^*(PAP) \to C^*_\env(PAP)$ which ensures that $\imath_E(x)^2 \geq (g^2 -\delta^2) P$ so that $x$ is a hermitian form with gap $\sqrt{g^2- \delta^2}$.

\item (Spectral truncations of quanta of geometry)
  As a special case of an invertible self-adjoint element we may also consider the operator $Y \in A$ describing the so-called quanta of geometry in \cite{CCM14}. It satisfies $Y^2 = 1$ so that $Y $ has spectrum contained in $\{ -1,1\}$. If $P$ is a projection, we find that the compressions $PYP$ define hermitian forms with gap $\sqrt{1-\delta^2}$ provided $\delta = \| [P,Y]\|<1$. 
\end{enumerate}
\end{ex}

  Let us continue the general treatment of hermitian forms in operator systems, and derive the following rigidity result. 
\begin{lma}
\label{lma:rigid}
Let $x$ be a hermitian form with gap $g>0$. If $y =y^* \in M_n(E)$ with  $\| x - y\| \leq \epsilon$ for some $\epsilon  < g^2/2 \|x\|$, then $y$ is a hermitian form with gap $\sqrt{g^2  - 2 \epsilon \| x\| }$. 
\end{lma}
\proof
The norm estimate implies that 
$$-2 \epsilon \| x\|  e \leq \imath_E (x)\imath_E (y-x) + \imath_E(y-x)\imath_E( x) \leq 2  \epsilon\| x\|e
$$
Hence by writing $y = x + (y-x)$ and noting that $\imath_E (y-x)^2 \geq 0$ we find that
\begin{align*}
  \imath_E(y)^2 &= \imath_E (x)^2 + \imath_E (x)\imath_E (y-x) + \imath_E(y-x)\imath_E( x) + \imath_E (y-x)^2
  \\ 
  & \geq ( g^2  - 2 \epsilon \| x\|  )e. \qedhere
\end{align*}
\endproof
We end this subsection by another standard operation of hermitian forms, which is their direct sum.
\begin{defn}
  \label{defn:direct-sum}
  Let $x \in H(E,n)$ and $x' \in H(E,n')$ be hermitian forms with respective gaps $g,g'$. Their {\em direct sum} is the hermitian form given by
  $$
x \oplus x' = \begin{pmatrix} x & 0 \\ 0 &  x' \end{pmatrix}\in H(E,n+n')
$$
which has gap equal to $\min\{ g,g'\}$. 
  \end{defn}

\subsection{Homotopy equivalence of hermitian forms}
One of the defining properties of an operator system is its matrix-order structure. This motivates us to consider homotopies of hermitian forms in $H(E,n)$ for a fixed matrix size $n \geq 0$, as follows:
\begin{defn}
  \label{defn:class-herm}
  Let $x,x' \in H(E,n)$ be hermitian forms. We say that $x \sim_n x' $ if there exists a hermitian form $\tilde x \in H(C([0,1]) \otimes E),n)$ such that
  $$
\tilde x(0) = x ; \qquad \tilde x(1) = x'
  $$
We denote the equivalence class of a hermitian form $x \in H(E,n)$ by $[x]_n$, and the set of equivalence classes of hermitian forms in $H(E,n)$ by $\V(E,n)$, or, equivalently,
$$
\V(E,n) = \pi_0 (H(E,n)).
$$
\end{defn}

This should be confronted with the usual notion of homotopy equivalence in, say, K-theory for $C^*$-algebras, where $n$ is allowed to vary (see also Section \ref{sect:semigroup} below). The following is immediate:
\begin{prop}
  \label{prop:ucp}
  Let $E$ and $F$ be unital operator systems. If $\phi:E \to F$ is a ucp map for which there exists a $*$-homomorphism $\tilde\phi :C^*_\env(E) \to C^*_\env(F)$ that makes the following diagram commute,
  \begin{equation}
    \label{eq:ucp-homom}
  \xymatrix {
    E \ar[r]^{\phi} \ar[d]_{\imath_E} & F  \ar[d]_{\imath_F}\\
    C^*_\env(E) \ar[r]^{\tilde \phi} & C^*_\env(F)
    }
  \end{equation}
  then there is an induced map $\phi^*:  \V(E,n) \to \V(F,n)$ defined by
  $$
\phi^* ([x]_n) = [\phi^{(n)}(x)]_n.
  $$
\end{prop}
\begin{corl}
  If $E$ and $F$ are completely order isomorphic then $\V(E,n) \cong \V(F,n)$ for all $n \geq 0$.
  \end{corl}

\begin{ex}
  \label{ex:V0-C}
  Consider $E=\C$. Then $\V(\C,n)$ is the set of homotopy equivalence classes of invertible hermitian $n \times n$ matrices. Any such matrix $x$ can be diagonalized with a unitary matrix, and since the unitary group $U(n)$ is connected, there is a homotopy between $x$ and the corresponding diagonal matrix. In turn, this diagonal matrix is homotopic (in the space of invertible matrices) to the corresponding signature matrix, which is unique up to ordering.  In other words, $\V(\C,n)$ can be parametrized by the signature $s$ of 
  the hermitian forms
  , yielding an isomorphism
  $$
  \V(\C,n) \cong \{  -n , -n+2 \ldots, n-2, n \}
  $$
  Note that the direct sum of two hermitian forms translates to the addition of the corresponding signatures: $(s,s')\mapsto s+s'$.
  
  Anticipating the discussion in Section \ref{sect:semigroup} there are maps $\imath_{nm} : \V(\C,n) \to \V(\C,m)$ for $n \leq m$ given by
  $$
\imath_{nm} ([x]_n) = [ x \oplus e_{m-n} ]_{m}.
  $$
In terms of the signature, 
we find that $\imath_{nm} (s) = s + m-n$. Moreover, there are commuting diagrams for any $m \geq n \geq 1$:
 \begin{equation*}
     \xymatrix {
       \V(\C,n) \ar[rd]_{\rho_n} \ar[rr]^{\imath_{nm}}&& \V(\C,m) \ar[ld]^{\rho_m}\\
&       \Z 
     }
  \end{equation*}
 where $\rho_n([x])$ is defined to be the so-called {\em negative index of inertia} of $x$, {\em i.e.} the number of negative eigenvalues of the hermitian form $x$. When expressed in terms of the signature $s$ of $x$, we have $\rho_n(s) =\frac 12 (n-s)$ from which commutativity of the diagram follows at once. This suggests that $\Z$ is the direct limit $\varinjlim \V(\C,n)$; we will come back to this soon (in Section \ref{sect:semigroup}).
  \end{ex}

\begin{lma}
  \label{lma:direct-sum}
Let $E$ and $F$ be unital operator systems. Then for any $n$ we have $\V(E\oplus F,n) \cong \V(E,n) \times \V(F,n)$.
\end{lma}
\proof
Let $(x,y)\in M_n(E) \oplus M_n(F)$ be a hermitian form; this has gap $g$ iff
$$
\imath_{E \oplus F}(x,y)^2 \geq g^2 \cdot (e_E, e_F) \iff \imath_E (x)^2 \geq g^2 e_E \text{ and }  \imath_E (y)^2 \geq g^2 e_F.
$$
Moreover, $(x,y) \sim (x',y')$ iff $x \sim x'$ and $y \sim y'$. 
From this the stated isomorphism follows. 
\endproof

The following result is a useful computational tool for $\V(E,n)$. Recall the multiplier $C^*$-algebra $A_E$ of the operator system $E$:
$$
A_E := \{ a \in C^*_\env(E): a E \subset E \} 
$$
\begin{lma}
  \label{lma:whitehead}
Let $u$ be a unitary in $M_n(A_E)$, and let $x \in M_n(E)$ be a hermitian form with gap $g$. Then $u x u^*$ is a hermitian form with gap $g$, and $u x u^* \oplus e_n$ is equivalent to $x \oplus e_n$ as hermitian forms in $H(E,2n)$.   
  \end{lma}
\proof
Clearly, $\imath_E(u x u^*)^2 = u \imath_E (x)^2 u^* \geq g^2 e$. For the homotopy equivalence, recall that by Whitehead's Lemma there exists a continuous path of unitaries $w_t \in M_{2n}(A_E)$ such that
$$
w_0 = \begin{pmatrix} 1 & 0 \\ 0 & 1  \end{pmatrix}; \qquad
w_1 = \begin{pmatrix} u & 0 \\ 0 & u^*  \end{pmatrix}.
$$
But then $w_t ( x \oplus e_n) w_t^*$ is a self-adjoint element in $C[0,1]\otimes M_{2n}(E)$, also with gap $g$, establishing an equivalence between $x \oplus e_n$ and $u x u^* \oplus e_n$, as desired. 
\endproof

\begin{prop}
  \label{prop:rigid}
  Let $x,y \in H(E,n)$ be hermitian forms both with gap $g$ for which $\| x - y\| < \epsilon$ such that $\epsilon < 4 g^2 $. Then there is a homotopy $\tilde x$ of hermitian forms in $C([0,1])\otimes E$ with gap $ \sqrt{g^2 - \epsilon^2/4}$ such that $\tilde x(0) = x$ and $\tilde x(1)=y$.
  \end{prop}
\proof
We claim that a homotopy of hermitian forms with the indicated gap is given by $\tilde x(t) = t x + (1-t)y$. It follows from $\| x - y\| < \epsilon$ that 
\begin{align*}
  & 
  (\imath_E(x-y))^2 \leq \epsilon^2 e,
\end{align*}
which implies that
$$
\imath_E(x) \imath_E(y)+\imath_E(y) \imath_E(x) \geq \imath_E(x)^2 + \imath_E(y)^2 - \epsilon^2 e.
$$
We will use this to find the gap of $\tilde x$ as follows:
\begin{align*}
 \imath_{C([0,1])\otimes E}( \tilde x)^2 &= t^2 \imath_E(  x)^2  + (1-t)^2   \imath_E(y)^2 + t (1-t) (\imath_E(x) \imath_E(y)+\imath_E(y) \imath_E(x))   \\
  & \geq (t^2 g^2 + (1-t)^2 g^2 +  t(1-t)(2 g^2 - \epsilon^2 ))e\\
  & \geq (g^2 -  \epsilon^2/4)e.\qedhere
  \end{align*}
\endproof

\begin{corl}
  \label{corl:rigid}
  Let $x \in H(E,n)$ be a hermitian form with gap $g$ and suppose that $y =y ^*\in M_n(E)$ is such that $\| x - y\| < \epsilon$ for some $\epsilon >0$ such that $g^2 > 2 \epsilon \|x \| + \epsilon^2/4$.
  Then there is a homotopy between $x$ and $y$ given by hermitian forms with gap $\sqrt{g^2- 2 \epsilon \| x \| - \epsilon^2/4}$.
\end{corl}
\proof
From Lemma \ref{lma:rigid} it follows that the element $y$ is a hermitian form with gap $\sqrt{g^2-2\epsilon \| x \|}$ and the same applies to $x$. By Proposition \ref{prop:rigid} it now follows that there is homotopy between $x$ and $y$ of hermitian forms with gap $\sqrt{g^2-2\epsilon \| x \|- \epsilon^2/4}$.
\endproof

\subsection{Semigroup of hermitian forms and $K_0$-group}
\label{sect:semigroup}
In the above, we have stressed the role of the matrix size $n \geq 0$ for the invariants $\V(E,n)$. However, in order to prepare for a comparison with $K$-theory for $C^*$-algebras we will need the limit object $\V(E) := \varinjlim \V(E,n)$ that we will now construct.

Given a unital operator system $(E,e)$ we consider the direct system of sets $(\V(E,n), \imath_{nm})$ where for $m \geq n$
\begin{align}
  \label{eq:dir-syst}
  \imath_{nm}: \V(E,n) &\to \V(E,m)\\
[x]_n & \mapsto [x \oplus e_{m-n}]_m, \nonumber
\end{align}
We denote the direct limit of the direct system \eqref{eq:dir-syst} by $\varinjlim \V(E,n)$. A more explicit description is given as follows: for $x \in H(E,n)$ and $x' \in H(E,n')$ we write $x \sim x'$ if there exists a $k \geq n,n'$ such that $ x \oplus e_{k-n} \sim_k x'\oplus e_{k-n'} $ in $H(E,k)$. We will write $[x]_E$, or simply $[x]$, for the equivalence class corresponding to $x \in H(E,n)$ and 
$\V(E) := \amalg_n \V(E,n)/_\sim$ for the corresponding set of equivalence classes. The following is then clear from the definition of the direct limit. 
\begin{prop}
  The set $\V(E)$ is the direct limit $\varinjlim \V(E,n)$ of the direct system \eqref{eq:dir-syst}. Moreover, it is a semigroup when equipped with the direct sum $[x]+ [x'] = [x \oplus x']$ and identity element $0 = [e]$. 
\end{prop}




We now arrive at our tentative definition of K-theory for unital operator systems.
\begin{defn}
  Let $(E,e)$ be a unital operator system. We define the $K$-theory group $K_0(E)$ of $E$ to be the Grothendieck group of $\V(E)$. 
\end{defn}

Before addressing the properties of the association $K_0(E)$ to an operator system, let us first check the consistency of our definition with the notion of $K_0$-groups when the unital operator system is in fact a unital $C^*$-algebra. 
\begin{prop}
  For a unital $C^*$-algebra $A$ the group $K_0(A)$ is isomorphic to the $C^*$-algebraic K-theory group of $A$. 
    \end{prop}
\proof 
Since a hermitian form $x$ in a $C^*$-algebra is invertible we may define a map between the semigroups (where $\mathcal P(A)$ denotes the semigroup of projections in $M_\infty(A)$ up to homotopy equivalence):
\begin{align*}
  \Phi: \V(A) &\to \mathcal P(A)\\
  [x] &\mapsto \left[p=\tfrac 12 (1-x |x|^{-1})  \right]
\end{align*}
This map is well-defined, since if $x \sim x'$ then also the corresponding projections $p,p'$ are homotopy equivalent. It also maps the neutral element $0 = [e]$ to $0=[0] \in \mathcal P(A)$, whilst being compatible with direct sums. Let us check it is bijective.

For injectivity, let $[x],[x'] \in \V(A)$ be mapped to $[p],[p']$ and assume that $[p]=[p']$. The homotopy $p_t$ of projections ($t \in [0,1]$) that implements the equivalence between $p_0 = p, p_1 = p'$ can be written as $p_t = \frac 12 (1-y_t)$ in terms of a self-adjoint element $y_t$ satisfying $y_t^2 = 1$, to wit $y_t = 1-2 p_t$. This a continuous family of hermitian forms such that $y_0 = x|x|^{-1}$ and $y_1 = x' |x'|^{-1}$  while also $x |x|^{-1} \sim x$ via the homotopy $x |x|^{-t}$ ($t \in [0,1]$) of hermitian forms, and the same applies to $x' |x'|^{-1} \sim x'$. Hence $x \sim x'$ which proves that $\Phi$ is injective.

Surjectivity follows by taking a $[p] \in \mathcal P(A)$ and defines a hermitian form $x = 1-2p$. Then $\Phi([x]) = [p]$ as desired. 
  \endproof
The next result prepares for invariance of $K_0$ under stable isomorphism (in Section \ref{sect:non-unit} below).
\begin{thm}[Stability of $K_0$]
  \label{thm:stab}
  Let $E$ be a unital operator system and let $N$ be a natural number. Then $\V(E)$ is isomorphic to $\V(M_N(E))$ (and so are the corresponding $K_0$-groups).
\end{thm}
\proof
We will give an explicit proof in the case $N=2$. For any $n>0$ we define a map $\jmath^{E,n} : M_n(E)\to M_n(M_2(E))$ by
\begin{align*}
  \jmath^{E,n} (x)= \begin{pmatrix} \begin{matrix} x_{11} & 0 \\ 0 & e \end{matrix}  & \rvline & \begin{matrix} x_{12} & 0 \\ 0 & 0 \end{matrix} & \rvline & \cdots & \rvline & \begin{matrix} x_{1n} & 0 \\ 0 & 0 \end{matrix}
    \\
    \hline
    \begin{matrix} x_{21} & 0 \\ 0 & 0 \end{matrix} & \rvline & \begin{matrix} x_{22} & 0 \\ 0 & e \end{matrix} & \rvline & \cdots & \rvline & \vdots
    \\
    \hline
    \vdots  & \rvline & \vdots & \rvline & \ddots & \rvline & \vdots
    \\
    \hline 
    \begin{matrix} x_{n1} & 0 \\ 0 & 0 \end{matrix}  & \rvline & \cdots  & \rvline & \cdots & \rvline & \begin{matrix} x_{nn} & 0 \\ 0 & e \end{matrix}
  \end{pmatrix}
      \end{align*}
Let $u_\sigma$ be the permutation matrix corresponding to the permutation
$$
\sigma = \begin{pmatrix} 1 & 2 & 3 & \cdots & n & n+1 & \cdots & 2n \\
  1 & 3 & 5 & \cdots &2n-1  & 2 & \cdots & 2n \end{pmatrix}
$$
This shuffles the columns and rows in such a way that 
$$
u_\sigma \cdot \jmath^{E,n}(x)\cdot u_\sigma^* = \begin{pmatrix} x & 0 \\ 0 & e_n \end{pmatrix}
$$
Since $u_\sigma \in M_{2n}(\C) \subseteq M_{2n}(A_E)$ it follows from Lemma \ref{lma:whitehead} that $\jmath^{E,n}(x) \sim x$, relative to $E$. Similarly, it follows that $\jmath^{E,n+n'}(x\oplus x') \sim \jmath^{E,n}(x) \oplus \jmath^{E,n'}(x')$, while the map $\jmath^{E,n}$ also respects non-degeneracy and self-adjointness. Note that also $\jmath^{E,n}(e_n) = e_{2n}$, and, in fact, $\jmath^{E,n}$ is unital completely positive. As such, it is contractive, hence continuous so that it maps homotopy equivalences to homotopy equivalences. We conclude from all of this that the induced map $\jmath^E: \V(E) \mapsto \V(M_2(E)), [x] \mapsto [\jmath^{E,n}(x)]$ is a well-defined morphism of semigroups.

Let us show that $\jmath^E$ is surjective: take $[y] \in \V(M_2(E),m)$. Then by the above it follows that $\jmath^{M_2(E),m}(y) \sim y$, relative to $M_2(E)$. Upon identifying $M_m(M_2(E))$ with $M_{2m}(E)$ there is an element $\tilde y \in H(E,2m)$ satisfying $\jmath^{E,2m}(\tilde y) \sim \jmath^{M_2(E),m} (y)$. But then $\jmath^E([\tilde y ]_E) = [y]_{M_2(E)}$ which shows surjectivity.

For injectivity  of $\jmath^E$ suppose $x \in H(E,n), x'\in H(E,n')$ are such that $\jmath^{E,n}(x) \sim  \jmath^{E,n'}(x')$, {\em i.e.} that $\jmath^E([x]) = \jmath^E([x'])$. Since $\jmath^{E,n}(x) \sim x$ and $\jmath^{E,n'}(x') \sim x'$ by the above, it follows that $x \sim x'$. But then $[x] = [x']$ which shows that $\jmath^E$ is injective. This completes the proof. 
\endproof
\begin{rem}
More generally, we may consider maps
\begin{align}
  \jmath^{E,n}_{NM} : M_n(M_N(E)) \to M_n(M_M(E)) \nn \\ 
  x \mapsto v_\sigma^* \cdot  \begin{pmatrix} x & 0 \\ 0 & e_{M-N} \end{pmatrix} \cdot v_\sigma,
  \label{eq:jMN}
\end{align}
in terms of a suitable permutation matrix $v_\sigma$ shuffling the rows and columns to identify $M_n(M_M(E))$ with $M_M(M_n(E))$. The corresponding maps $\jmath_{NM}: \V(M_N(E)) \to \V(M_M(E))$ then yield a direct system, and the above Theorem generalizes to an isomorphism $\varinjlim \V(M_N(E)) \cong \V(E)$. 
\end{rem}

\subsection{$K_0$ and maps between operator systems}
The behavior of hermitian forms with respect to ucp maps as obtained in Proposition \ref{prop:ucp} translates into the following: 
\begin{prop}
\label{prop:K0-functor}
Let $E,F$ be unital operator systems and $\phi:E \to F$ a restriction of a $*$-homomorphism between the corresponding $C^*$-envelopes (so that Equation \eqref{eq:ucp-homom} is satisfied). Then the induced map $\phi^* : \V(E) \to \V(F)$ of semigroups is well-defined and induces a map between the corresponding $K_0$-groups (also denoted $\phi^*$).
\end{prop}
\proof
We already know from Proposition \ref{prop:ucp} that hermitian forms are mapped to hermitian forms. Clearly the map $\phi$ behaves well with respect to direct sums of hermitian forms, and unitality of $\phi$ implies that $[e_E]$ is mapped to $[e_F]$. Finally, as in Proposition \ref{prop:ucp} homotopies are mapped to homotopies.
\endproof

\begin{corl}
  \label{corl:K0-coi}
Let $E,F$ be unital operator systems. If $E$ and $F$ are completely order isomorphic, then $K_0(E) \cong K_0(F)$.
  \end{corl}
\proof
This follows from the previous Proposition in combination with the fact that a unital complete order isomorphism $\phi: E \to F$ extends to a unital $*$-isomorphism $\tilde  \phi: C^*_\env(E) \to C^*_\env(F)$ \cite[Theorem 4.1]{Ham79}. 
\endproof

The following is then immediate. 
\begin{prop}
  \noindent 
  \begin{enumerate}
  \item For every unital operator sytem $E$, $\id_E^* = \id_{K_0(E)}$;
  \item If $E, F,G$ are unital operator systems, and if $\phi:E \to F, \psi:F \to G$ are restrictions of $*$-homomorphisms between the respective $C^*$-envelopes, then $\psi^* \circ \phi^*= \psi^* \circ \phi^*$;
  \item $K_0(\{0\}) = \{ 0\}$;
    \end{enumerate}  
  \end{prop}

We also record the following behaviour of $K_0$ with respect to direct sums of operator systems, which is a direct consequence of Lemma \ref{lma:direct-sum}.
\begin{prop}
  \label{prop:direct-sum-K0}
  Let $E,F$ be unital operator systems. Then $K_0(E \oplus F) \cong K_0(E) \times K_0(F)$. 
  \end{prop}

\section{Non-unital operator systems and stability of K-theory}
\label{sect:non-unit}

\subsection{Non-unital operator systems}
Recall \cite{Wer02} ({\em cf.} \cite{CS20}) that the partial unitization of a non-unital operator system $E$ is given by the $*$-vector space $E^\sharp = E \oplus \C$ with matrix order structure:
$$
(x,A) \geq 0 \text{ iff } A \geq 0  \text{ and } \phi(A_\eps^{-1/2} x A_\eps^{-1/2} ) \geq -1
$$
for all $\eps >0$ and noncommutative states $\phi \in \mathcal S_n(E)$, and where $A_\eps = \eps \mathbb I_n + A$. The matrix order units $e^\sharp_n$ in $M_n(E^\sharp)$ are given by the identity matrices $\mathbb I_n \in M_n(\C)$. This turns $(E^\sharp, e^\sharp)$ into a unital operator system. Moreover, given a completely contractive completely positive map $\phi:E \to F$ we have that the canonical extension $\phi^\sharp: E^\sharp \to F^\sharp$ is a ucp map (\cite[Lemma 4.9]{Wer02}) . 

We now extend the definition of K-theory to non-unital operator systems as follows: 

\begin{defn}  Let $E$ be a non-unital operator systems and $E^\sharp$ its unitization. We define the sets 
 \begin{align*}
    \tilde \V(E,n) &:= \pi_0 \left( \left \{ (x,A) \in H(E^\sharp,n) : A \sim_n \mathbb I_n \right\} \right)
 \end{align*}
 The Grothendieck group of the direct limit semigroup is denoted by $\tilde K_0(E)$.
\end{defn}

\begin{prop}
Let $E$ and $F$ be operator systems. If $E$ and $F$ are completely isometric, completely order isomorphic, then $\tilde \V(E) \cong \tilde \V(F)$ as semigroups. Consequently, in this case $\tilde K_0(E) \cong \tilde K_0(F)$. 
\end{prop}
\proof
This follows from the fact that for a completely isometric, complete order isomomorphism $\phi:E \to F$ the induced map $\phi^\sharp :E^\sharp \to F^\sharp$ is a complete order isomorphism which furthermore respects the property $A \sim_n \mathbb I_n$. 
\endproof

Let us also compare this with with our previous definition of K-theory in the case of {\em unital} operator systems. 

\begin{prop}
  \label{prop:K0-unit}
  Let $E$ be a unital operator system and let $E^\sharp$ be its partial unitization. Then for any $n \geq 1$ there are isomorphisms $\tilde \V(E,n) \cong \V(E,n)$. Consequently, in this case $\tilde K_0(E) \cong K_0(E)$.
  \end{prop}
\proof
In the unital case, there is a unital complete order isomorphism between $E^\sharp$ and $E \oplus \C$ given by $(x, \lambda) \mapsto (x+\lambda e, \lambda)$, where we have equipped $E \oplus \C$ with the induced direct sum order structure \cite[Lemma 4.9(b)]{Wer02}. By Lemma \ref{lma:direct-sum} we have that $H(E \oplus\C,n) = H(E,n) \times H(\C,n)$ from which it follows that
$$
\left \{ (x,A) \in H(E^\sharp,n) : A \sim_n \mathbb I_n \right\}  \cong H(E,n) \times \left \{ A \in H(\C,n) : A \sim_n \mathbb I_n \right\}  
$$
This isomorphism clearly respects direct sums and the unit, so that taking homotopy equivalence classes of this set of hermitian forms yields the statement. 
\endproof

\subsection{Stability of $K_0$}
One of the crucial features of $K$-theory for $C^*$-algebras is that it is Morita invariant, or, equivalently, invariant under stable isomorphism. We will now establish that $K$-theory for unital operator systems shares this property, {\em i.e.} we will relate the $K_0$-groups of $E$ and $\mathcal K \otimes E$. 

Consider a direct system of Hilbert subspaces $\{ P_N \H \}_{N \geq 0}$ in $\H$ with $P_N \H \cong \C^N$ so that we realize $\K = \K(\H) \cong \varinjlim M_N(\C)$. We then obtain the stabilization of $E$ by the following series of completely contractive, completely positive maps defined for $M \geq N$:
$$
\kappa_{NM} : M_N(E) \to M_{M}(E), \qquad x \mapsto \begin{pmatrix} x& 0 \\ 0 &  0_{M-N} \end{pmatrix},
$$
The inductive limit of the sequence $(M_N(E),\kappa_{NM})$ is $\mathcal K \otimes E$ with connecting maps $\kappa_{N,\infty} : M_N(E) \to \mathcal K \otimes E$.

\begin{lma}
  The maps 
  $(\kappa_{NM}^\sharp)_*: \tilde \V (M_N(E),n) \to \tilde \V(M_M(E),n)$ induced by $\kappa_{NM}$ make the following diagram commute for any $n$:
$$
\xymatrix{
  \tilde \V (M_N(E),n) \ar[d]_{\cong}\ar[rr]^-{ (\kappa_{NM}^\sharp)_*} && \tilde \V (M_M(E),n)\ar[d]_{\cong}\\
   \V (M_N(E),n) \ar[rr]^-{\jmath_{NM} } && \V (M_M(E),n)\\
}
  $$
where $\jmath_{NM}$ is defined in Equation \eqref{eq:jMN} and the vertical isomorphisms are the ones from Proposition \ref{prop:K0-unit}.
\end{lma}
\proof
The maps $\kappa_{NM}^\sharp: M_N(E)^\sharp \to M_M(E)^\sharp$ and their amplifications are given explicitly by
$$
(\kappa_{NM}^\sharp)^{(n)} : M_n(E^\sharp) \to M_n( M_{N}(E)^\sharp) , \qquad (x,A) \mapsto  \left( u_\sigma \begin{pmatrix} x& 0 \\ 0 &  0_{n(M-N)}  \end{pmatrix} u_\sigma^* , A \right).
$$
Here $u_\sigma$ is a permutation matrix similar to the one appearing in the proof of Theorem \ref{thm:stab}: it  identifies $M_M(M_n(E))$ with $M_n(M_M(E))$ by a suitable shuffle of rows and columns. 
When we identify $M_N(E)^\sharp$ with the direct sum operator systems $M_N(E) \oplus \C$ ({\em cf.} Proposition \ref{prop:K0-unit}), then the map $(\kappa_{1N}^\sharp)^{(n)}$ becomes
$$
(x,A) \in M_n(M_N( E)) \oplus M_n(\C) \mapsto \left(u_\sigma\left( \begin{smallmatrix} x & 0 \\ 0 & A e_{M-N} \end{smallmatrix}\right) u_\sigma^*, A \right) \in M_n( M_{M}(E)) \oplus M_n(\C) 
$$
If we assume that $A \sim_n \I_n$ we see that up to homotopy this map coincides with $\imath_{NM}$, which completes the proof.
\endproof

\begin{thm}
For a unital operator system $E$ we have $\tilde K_0(\mathcal K \otimes E) \cong K_0(E)$. 
  \end{thm}
\proof
The above Lemma, in combination with Theorem \ref{thm:stab} and Proposition \ref{prop:K0-unit}, yields that $\varinjlim \tilde \V (M_N(E)) \cong \V(E)$. We will show that the universal map $u: \varinjlim \tilde \V(M_N(E)) \to \tilde \V(\K \otimes E)$ for the direct limit in the following diagram is an isomorphism:
$$
\xymatrix{
  \tilde \V (M_N(E)) \ar[rd] \ar@/^-2.0pc/[rdd]_{(\kappa_{N\infty}^\sharp)_*}  \ar[rr]^{(\kappa_{NM}^\sharp)_*} && \tilde \V (M_M(E)) \ar[ld] \ar@/^2.0pc/[ldd]^{(\kappa_{M\infty}^\sharp)_*} \\
& \varinjlim \tilde \V (M_N(E)) \ar@{-->}[d]_{u}&  \\
& \tilde \V(\K \otimes E)
}
$$

For injectivity of $u$, take $[(x,A)] \in \tilde \V(M_N(E),n)$ and $[(x',A')] \in \tilde \V(M_N(E),n')$ for some $N$ and $n,n'$. 
Assume that $[(\kappa_{N,\infty}^{(n)}(x),A)] = [(\kappa_{N,\infty}^{(n')}(x),A)] $ as elements in $\tilde \V (\K\otimes E)$. In other words, there exists a family $(\tilde x,\tilde A) \in H ( C[0,1] \otimes (\mathcal K \otimes E)^\sharp ,k)$ for some $k \geq n,n'$ such that the following hold:
\begin{gather*}
  \tilde x(0) = \kappa_{N,\infty} ^{(n)}(x) \oplus 0_{k-n} , \qquad\tilde A(0)=A \oplus \I_{k-n}, \\
  \tilde x(1) = \kappa_{N,\infty}^{(n)} (x') \oplus 0_{k-n'}  \qquad \tilde A(1) = A' \oplus \I_{k-n'}.
  \end{gather*}
Consider the compressions $R_M: \mathcal K \otimes E  \to M_M(E)$ given by $y \mapsto (P_M \otimes e)y (P_M \otimes e)$. Then for the above $x \in M_n(M_N(E))$ we have
$$
\kappa_{NM}^{(n)}(x) = R_M^{(n)} (\kappa_{N \infty}^{(n)} (x)) \equiv  (P_M \otimes e_n) \kappa_{N \infty}^{(k)}(x) (P_M \otimes e_n).
$$
and similarly for $x'$. We may also compress the homotopy to give $(R_M^\sharp)^{(k)}((\tilde x, \tilde A))  = (R_M^{(k)} (\tilde x),\tilde A )$ with end points:
\begin{gather*}
  (R_M^\sharp)^{(k)}( ( (\tilde x(0)),\tilde A (0)) = \left( \kappa_{NM}^{(n)}(x) \oplus 0_{k-n},  A \oplus \I_{k-n} \right)\\
  (R_M^\sharp)^{(k)}( ( (\tilde x(1)),\tilde A (1)) = \left( \kappa_{NM}^{(n')}(x') \oplus 0_{k-n'},  A' \oplus \I_{k-n'} \right)
\end{gather*}
A priori the family $(R_M^\sharp)^{(k)}((\tilde x, \tilde A))$ lies in $M_k(C[0,1] \otimes  (M_M( E))^\sharp)$ so we need to show that it completely lies in $\tilde \V( M_M(E),n)$, at least for some sufficiently large $M$.
In other words, we need to show that the gap of $(R_M^{(k)} (\tilde x),\tilde A)$  is strictly positive, while $\tilde A(t) \sim \mathbb I_k$ for all $t$. The latter fact is clear from the assumption that both $A \sim \I_n$ and $A' \sim \I_{n'}$. For the former claim,  we identify $C[0,1] \otimes (M_M (E))^\sharp$ with $C[0,1] \otimes (M_M (E) \oplus \C)$ which maps $(R_M^{(k)} (\tilde x),\tilde A  )$ to $(R_M^{(k)} (\tilde x) + P_M \otimes \tilde A e,\tilde A)$ since $P_M \otimes e$ is the order unit of $M_M(E)$. In order to check that the latter has strictly positive gap, we realize $M_k((\K \otimes E)^\sharp)$  as concrete operators in $B( \H \otimes ( \H')^{\oplus k})$ for some Hilbert space $\H'$. We compute
\begin{align*}
  (R_M^{(k)}(\tilde x )+ P_M \otimes \tilde A  )^2 &= ((P_M \otimes e_k) (\tilde x + \mathbb I_\H \otimes \tilde A )(P_M \otimes e_k)   )^2\\
  & =  ((P_M \otimes e_k) (\tilde x + \mathbb I_\H \otimes \tilde A )^2(P_M \otimes e_k)   \\
  &\qquad +  (P_M \otimes e_k  )[P_M \otimes e_k ,\tilde x + P_M \otimes \tilde A ]^2 (P_M \otimes e  )\\
  & \geq g^2(P_M \otimes e_k)  -\| [P_M \otimes e_k , \tilde x]\| ^2.
\end{align*}
We have used that $(\tilde x(t),\tilde A(t)) \in H( (\K \otimes E)^\sharp,k)$ has gap given by some $g>0$ for all $t$. 
Moreover, $\| [P_M \otimes e_k , \tilde x ]\| \to 0 $ as $N \to \infty$ so that it follows that $(R_M^{(k)}(\tilde x),\tilde A )$ is a hermitian form. We thus find that $[(R_M^{(k)}(\tilde x)(t),\tilde A(t) )] \in \tilde \V( M_M(E),k)$ for all $t$, so that $[(\kappa_{NM}^{(n)}(x),A)] = [(\kappa_{NM}^{(n')}(x'), A)]$ in $\tilde \V(M_M(E))$. Since we also know that $\kappa_{NM}^\sharp$ induces isomorphisms between $\tilde \V(M_N(E))$ and $\tilde \V (M_M(E))$ we conclude that $[(x,A)] = [(x', A)]$ as elements in $\tilde \V(M_N(E))$.

For surjectivity, take an arbitrary $[(x,A)] \in \tilde \V( \K \otimes E, n)$. We may approximate $x$ by a finite rank operator (for instance using the compression $R_N^{(n)}$) so that for all $\epsilon>0$ there exists a $x_0\in M_n(M_N(E))$ for some $N$ so that
$$
\| x - \kappa_{N,\infty}^{(n)}(x_0) \| < \epsilon.
$$
We then also have
$$
\| (x,A) -  (\kappa_{N,\infty}^\sharp)^{(n)} (x_0,A) \|< \epsilon
$$
so that if $(x,A)$ has gap $g$ and we choose $\epsilon$ small enough so that $g^2 - \epsilon \| x\| -  \epsilon^2/4 >0$, then $(x,A)$ and $(\kappa_{N,\infty}^\sharp)^{(n)} (x_0,A)$ are homotopy equivalent by Corollary \ref{corl:rigid}. The element $(x_0,A)$ is the sought-for element in $\tilde \V(M_N(E)) \cong \tilde \V(E)$ that maps to $[(x,A)]$.
\endproof

\begin{corl}
Let $E$ and $F$ be Morita equivalent unital operator systems. Then $K_0(E) \cong K_0(F)$. 
  \end{corl}
\proof
In \cite{EKT21} it is shown that $E$ and $F$ are Morita equivalent whenever $\K \otimes E \cong \K \otimes F$ via a completely isometric, complete order isomorphism.
\endproof

\section{Application to the spectral localizer}
\label{sect:appl}
In \cite{LS18a,LS18b} the spectral localizer was introduced as a powerful tool for computing index pairings. They relate a certain Fredholm index to the signature of a finite-dimensional matrix ---the so-called spectral localizer. We will put it in the context of our notion of K-theory, realizing the spectral localizer as a 
map $\V(E,n) \to \Z$. 
the spectral localizer, including its relation to spectral flow, we refer to the excellent textbook \cite{DSW23} and references therein.

In order to describe the index map on $\V(E,n)$, we start with an operator system spectral triple \cite{CS20}. 
\begin{defn}
  A (unital) {\em operator system spectral triple} is given by a triple $(E,\H,D)$ where $E$ is a unital operator system realized concretely so that $E \subseteq C^*_\env(E) \subseteq B(\H)$, 
  and a self-adjoint operator $D: \dom(D) \to \H$ such that 
  \begin{itemize}
  \item the commutators $[D,x]$ extend to bounded operators for all $x \in \E$ for a dense $*$-subspace $\E \subseteq E$;
    \item the resolvent $(i+D)^{-1}$ is a compact operator.
  \end{itemize}
  An operator system spectral triple is called {\em even} if in addition to the above, there is a grading operator $\gamma$ on $\H$ (so that $\gamma^* = \gamma, \gamma^2 = 1_\H$) which commutes with all $x \in E$ and anti-commutes with $D$. Otherwise, it is called {\em odd}.
\end{defn}
In the even case, we can decompose $\H = \H_+ \oplus \H_-$ according to the eigenvalues of $\gamma$, and decompose accordingly
$$
D = \begin{pmatrix} 0 & D_0 \\ D_0^* & 0 \end{pmatrix}.
$$

Given an even operator system spectral triple, a parameter $\kappa >0$ and a hermitian form $x \in H(\E,n)$ we now define the {\em even spectral localizer} \cite{LS18a} as
\begin{align}
  L_{\kappa} (D,x) &= \begin{pmatrix} x & \kappa (D_0)^{\oplus n}  \\ \kappa (D_0^*)^{\oplus n} & - x  \end{pmatrix}.
\end{align}


\begin{prop}
  \label{prop:homotopy-sf}
  Let $(E,\H,D)$ be an even finite-dimensional operator system spectral triple. 
  Let $x \in H(\E,n)$ with gap $g>0$ and let $\kappa_0 = g^2  \| [D,x] \|^{-1}$. Then 
  the index map defined by the signature of the spectral localizer, 
     $$
\ind_D([x]) = \frac 12 \Sig (L_\kappa(D,x)),
       $$
is constant for all $\kappa < \kappa_0$ and invariant under homotopy equivalence. Consequently, it induces a map $\ind_D: \V(\E,n) \to \Z$.
  \end{prop}
\proof
Without loss of generality we take $n=1$ and compute very similar to \cite{LS18a} that
\begin{equation}
L_\kappa(D,x)^2 = \begin{pmatrix} \kappa^2 D_0 D_0^* + x^ 2 & \kappa [D_0,x] \\ \kappa [D_0,x]^* & \kappa^2 D_0^* D_0 +x^2 \end{pmatrix} \geq \left( g^2 -\kappa \| [D,x]\| \right) 1_{M_2(B(\H))}.
\label{eq:loc-square}
\end{equation}
Hence $L_\kappa(D,x)$ is invertible (and thus has well-defined signature) provided $\kappa< \kappa_0$ and, moreover, the signature is constant for all $\kappa <\kappa_0$ as no eigenvalues will cross the origin.

Also, one may easily check that
$$
\Sig (L_\kappa(D,x \oplus x')) = \Sig (L_\kappa(D,x ))+ \Sig (L_\kappa(D,x' )).
$$
In order to see that $\Sig (L_\kappa(D,e )) = 0$ consider an eigensystem $\{v_\lambda =  (v_\lambda^+, v^-_\lambda)\}_\lambda$ in $\H_+ \oplus \H_-$ such that $D_0 v^-_\lambda = \lambda v^+_\lambda$ and $\gamma v_\lambda^\pm = \pm v_\lambda^\pm$. We may write the matrix of $L_\kappa$ in this eigensystem as
$$
\langle v_\lambda, L_\kappa(D,x) v_{\lambda'} \rangle  =\delta_{\lambda\lambda'}\begin{pmatrix} 1 & \lambda \\ \lambda & -1 \end{pmatrix} \sim \delta_{\lambda\lambda'} \begin{pmatrix} \sqrt{1+\lambda^2} & 0 \\ 0& -\sqrt{1+\lambda^2} \end{pmatrix} 
$$
and the signature of the latter matrix vanishes.

Consider now a homotopy $\tilde x$ in $H(E,n)$ between $x = \tilde x(0)$ and $x' = \tilde x (1)$. We define $\kappa_0' = g^2 / \sup_t \| [D,\tilde x(t)]\|$ and note that for all $\kappa < \kappa_0'$ 
$$
L_{\kappa}(D,\tilde x)^2 \geq g^2 - \kappa \sup_t \| [D, \tilde x(t)]\| 
$$
by a computation similar to the one in Eq. \eqref{eq:loc-square}. So as long as $\kappa < \kappa_0'$ we find that $\Sig L_{\kappa}(D,\tilde x(t))$ is constant in $t$. We combine this with the fact that $\kappa < \kappa_0' < \kappa_0$ to conclude that $\Sig L_{\kappa}(D,\tilde x(0)) = \Sig L_{\kappa}(D,\tilde x(1))$. This completes the proof. 
\endproof

We may now rephrase the main results of \cite{LS18a}. For an invertible self-adjoint element $x \in A$ one considers the class $[p=\frac 12 (1-x|x|^{-1})]$ in the K-theory of the $C^*$-algebra $A$. Given a spectral triple $(A,\H,D)$ one would like to compute the index $\ind pDp$. As shown in {\em loc.cit.} this index may be computed in terms of a spectrally truncated (operator system) spectral triple $(P A P,P \H , P D P)$ for a spectral projection $P$ (of $D$) of sufficiently high rank. Indeed, we then have
\begin{equation}
\ind (pDp)= \ind_{PDP}  ([P x P ])
     \label{eq:ind-loc}
\end{equation}
where the map on the right-hand side is the index map on $\V(PAP,n)$ that appeared in Proposition \ref{prop:homotopy-sf}. 

\subsection{Example: spectral localizer on the torus}
\begin{figure}
\includegraphics[scale=.2]{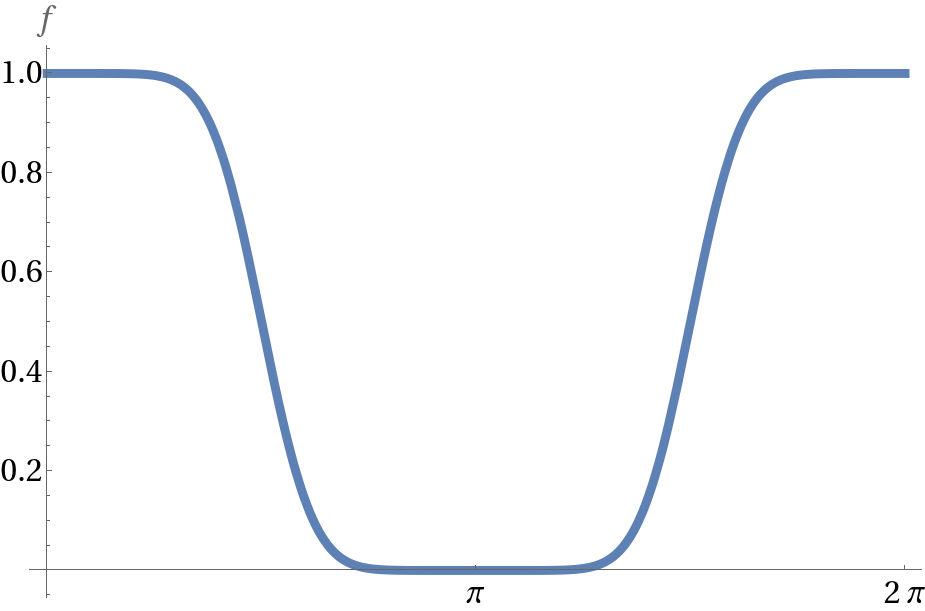}
\\
\includegraphics[scale=.2]{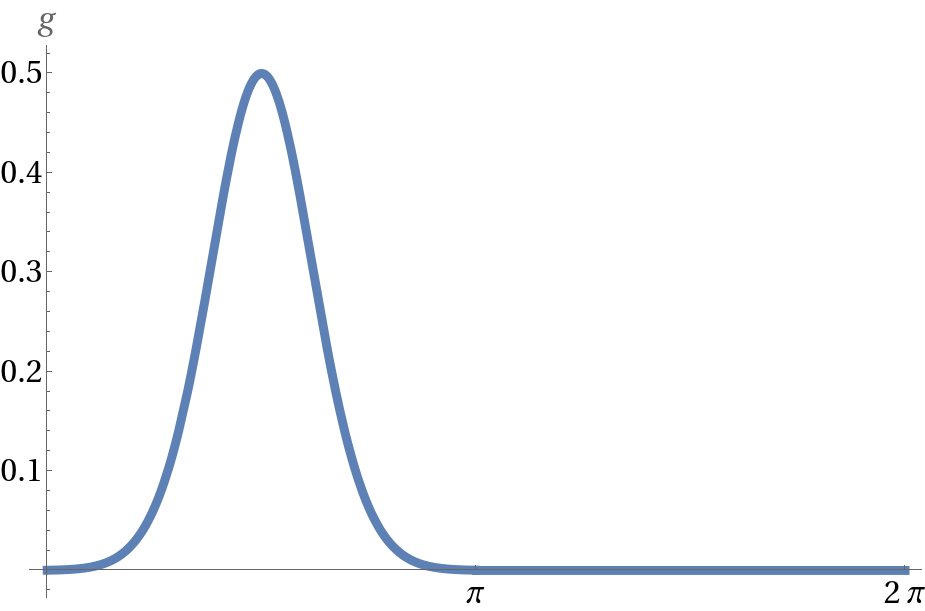}\hspace{1cm}
\includegraphics[scale=.2]{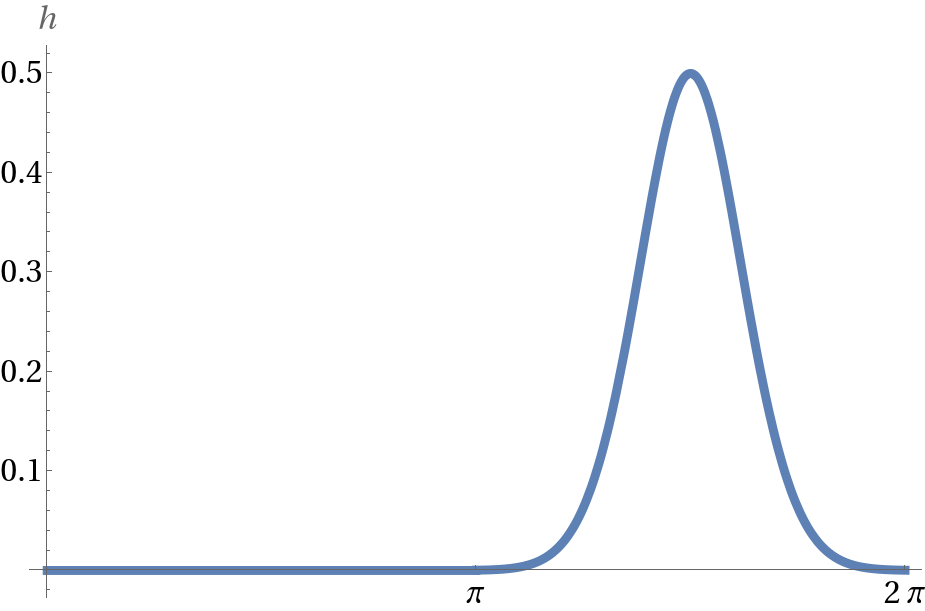}
\caption{Functions $f,g,h$ that ensure that $p$ in \eqref{eq:proj-T2} is a projection.}
\label{fig:proj-T2}
\end{figure}

Let us illustrate the spectral localizer for a spectrally truncated two-torus. First, recall from \cite{Lor86} the class of projections on the torus ({\em cf.} \cite[Chapter 6]{Sui14}, much inspired by the so-called Powers--Rieffel projections on the noncommutative torus \cite{Rie81}:
\begin{equation}
\label{eq:proj-T2}
p = \begin{pmatrix} f & g+hU^* \\ g+ hU & 1-f \end{pmatrix} \in M_2(C(\T^2)),
\end{equation}
where $f,g,h$ are real-valued (periodic) functions of the first variable $t_1$, and $U$ is a unitary depending only on the second variable $t_2$, say $U(t_2)=e^{i m t_2}$. The projection property $p^2=p$ translates into the two conditions
$$
gh =  0, \qquad g^2 + h^2 = f-f^2. 
$$
A possible solution of these relations is given by
$$
0 \leq f \leq 1 \quad \text{such that } f(0)=1, \quad f(\pi) = 0 
, 
$$
and then $g = \chi_{[0,\pi]} \sqrt{f-f^2}$ and $h = \chi_{[\pi,2\pi]} \sqrt{f-f^2}$, where $\chi_X$ is the indicator function for the set $X$ (see Figure \ref{fig:proj-T2}).
The Dirac operator on the two-torus is defined on the core $C^\infty(\T^2) \otimes \C^2$ by
$$
D= \begin{pmatrix} 0 & D_0 \\ D_0^* & 0 \end{pmatrix}; \qquad D_0 = i \partial_{t_1} + \partial_{t_2}.
$$
The corresponding spectrum is $\left\{ \pm \sqrt{n_1^2+n_2^2}\right\}_{n_1,n_2\in \Z}$.

Let us now consider a spectral projection corresponding to a (discrete) ball of radius $\rho$, {\em i.e.} $P_\rho = \chi(|D|\leq \rho)$. We  consider the corresponding spectral compression of $Y=1-2p$, {\em i.e.}
$$
P_\rho  Y P_\rho  = \begin{pmatrix} P_\rho -2 P_\rho  fP_\rho  & -2P_\rho gP_\rho -2P_\rho h U P_\rho   \\ -2P_\rho gP_\rho -2P_\rho h U^*P_\rho  & -P_\rho  +2P_\rho fP_\rho  \end{pmatrix} \in M_2(P_\rho  C^\infty(\T^2)P_\rho )
  $$
For suitable $P_\rho $ these are hermitian forms in the truncated operator system, that is, $P_\rho YP_\rho  \in \V(P_\rho  C(\T^2)P_\rho ,2)$.

\begin{figure}
  \includegraphics[scale=.38]{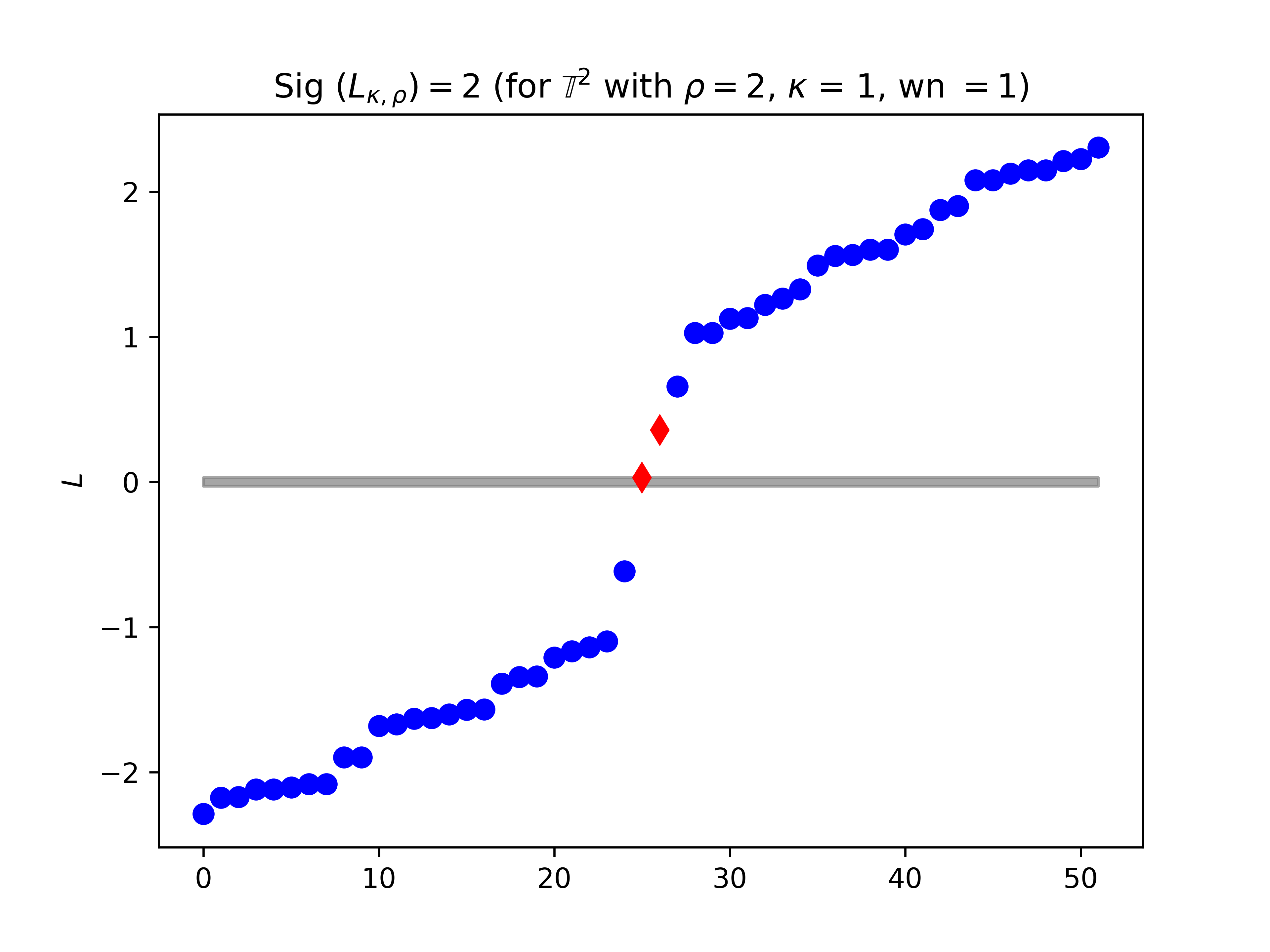}
  \includegraphics[scale=.38]{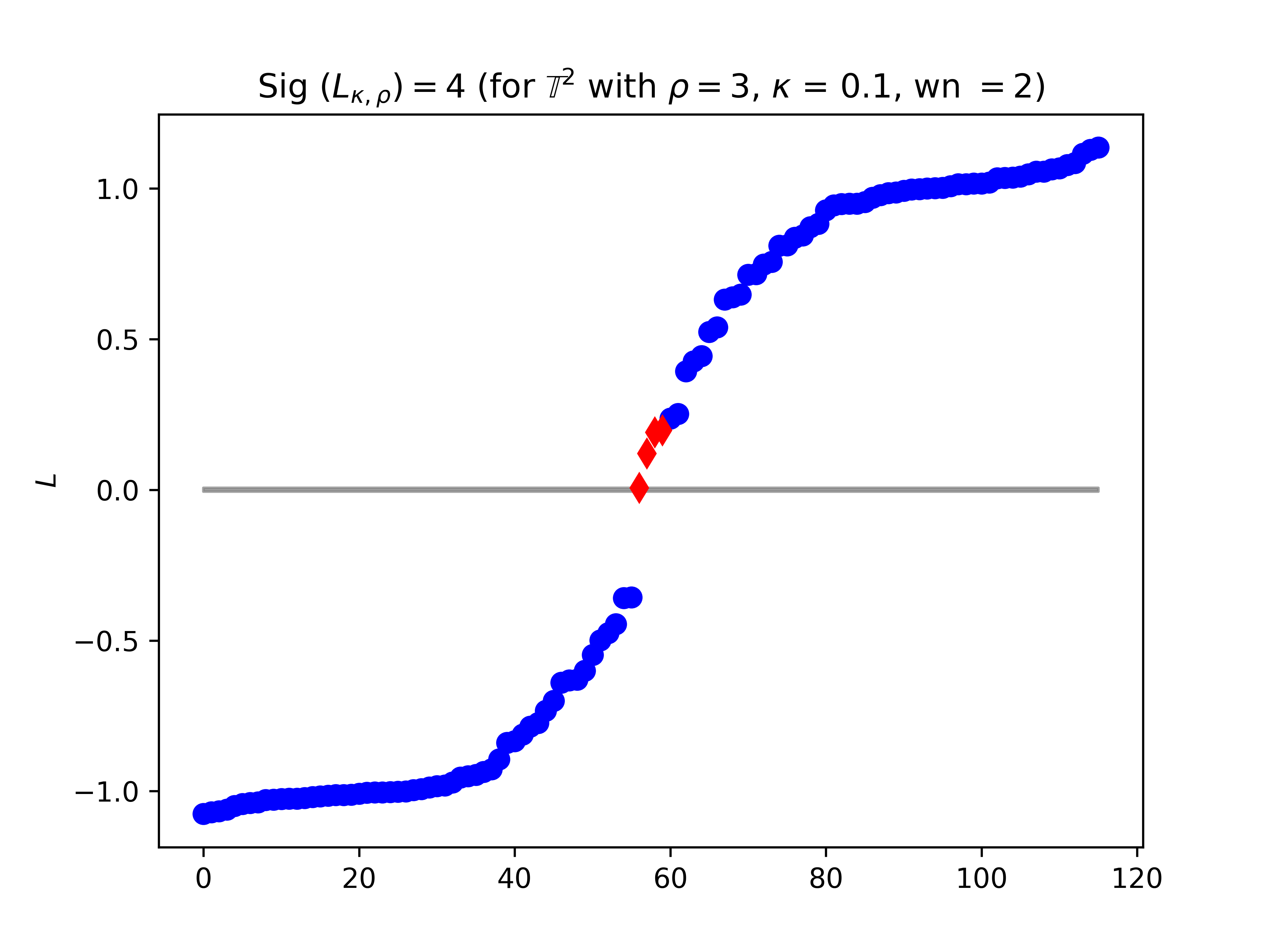}
  \caption{The negative and positive eigenvalues for the spectral localizer $L_{\kappa,\rho}$ on the torus for $m=1, \kappa= 1, \rho=2$ (left) and  $m=2, \kappa= 0.1, \rho=3$ (right). The red (diamond-shaped) dots indicate the surplus of positive eigenvalues as compared to the negative ones. }
    \label{fig:L-T2}
\end{figure}
The spectral localizer is now given by the following matrix: 
$$
L_{\kappa,\rho}:=
L_{\kappa}(P_\rho D P_\rho,P_\rho Y P_\rho)  = \begin{pmatrix} P_\rho YP_\rho  & \kappa P_\rho  D_0 P_\rho  \\ \kappa P_\rho   D_0^* P_\rho  &  -P_\rho YP_\rho  \end{pmatrix}
  $$

The main result of \cite{LS18a} ({\em cf.} Eq. \eqref{eq:ind-loc}) now implies that for suitable $\rho,\kappa$ the index of $pD p$ can be expressed in terms of the signature of the above spectral localizer.  More precisely, in terms of the above index map on $\V(P_\rho C(\T^2)P_\rho , 2)$ we have:
$$
\ind (pDp) = \ind_{P_\rho DP_\rho }([P_\rho YP_\rho ])  \equiv \frac 12 \text{Sig} ~L_{\kappa,\rho}
$$
In Figure \ref{fig:L-T2} we illustrate the resulting signature by showing the negative and positive eigenvalues of $L_{\kappa,\rho}$. We find already for low $\rho$ that the signature of the spectral localizer is equal to (twice) the winding number $m$, that is to say, the index of $pD p$.

\section{Outlook}
We have proposed a generalization of K-theory for $C^*$-algebras to operator systems, which in the spirit of Witt was based on hermitian forms. Several aspects are still to be developed, and which we leave for future research. This includes a careful study of the functorial properties of $K_0$ and, in particular, how $K_0$ behaves with respect to approximations results. For instance for approximations of $C^*$-algebra by finite-dimensional algebras in \cite{BK97,CW23} where the maps are approximately multiplicative, or approximately order-zero cpc maps. Or in the context of quantum metric spaces where it turns out that small Hausdorff distance between two such spaces ---the state spaces of two $C^*$-algebras $A,B$--- allows one to relate projections in $A$ to projections in $B$ \cite{Rie10,Rie18}. 

Another open problem is to identify the higher K-groups, starting with $K_1$. Again, one needs to find a analogue of the assumption of being unitary. However, just demanding invertibility appears too weak; one can easily show that the set of all invertible elements in a finite-dimensional operator system is contractible ---arguing much as in the case of $GL(n,\C)$. Nevertheless, from the odd spectral localizer \cite{LS18b} we know that some index-theoretic information is maintained after spectrally truncating a unitary. Again this fact will be leading in the development of $K_1$, which we also leave for future research. This extends of course to finding an analogue of Bott periodicity. 

Dually, one may be interested also in defining K-homology for operator systems. Note here that the notion of a Fredholm module makes perfect sense for unital operator systems realized concretely in Hilbert space, but that this immediately gives rise to a Fredholm module for the pertinent $C^*$-extension. The task is thus to find the right notion of equivalence.

\newcommand{\noopsort}[1]{}\def\cprime{$'$}

\end{document}